\documentclass[twoside,11pt]{article}

\usepackage{jmlr2e}

\usepackage[utf8]{inputenc} 
\usepackage[T1]{fontenc}    
\usepackage{hyperref}       
\usepackage{url}            
\usepackage{booktabs}       
\usepackage{amsfonts}       
\usepackage{nicefrac}       
\usepackage{microtype}      
\usepackage{times}
\usepackage{graphicx} 
\usepackage{subcaption}
\usepackage{tabularx}
\usepackage{bbm}
\usepackage{url}
\usepackage{amsmath}
\usepackage{bbm}
\usepackage{caption}
\usepackage{float}
\usepackage{amssymb}
\usepackage{grffile}
\usepackage{epstopdf}
\usepackage[colorinlistoftodos]{todonotes}
\usepackage{multirow}

\usepackage{algorithm}
\usepackage{algorithmic}
\usepackage[algo2e,linesnumbered,lined,boxed,vlined]{algorithm2e}

\usepackage{color}

\usepackage[normalem]{ulem}

\newcommand{\hide}[1]{}

\usepackage[nobiblatex]{xurl}

\usepackage[numbers]{natbib}

\begin{document}

\title{On Philomatics and Psychomatics \\for Combining Philosophy and Psychology with Mathematics
}



\author{\name Benyamin Ghojogh \email bghojogh@uwaterloo.ca 
\\
\name Morteza Babaie \email mbabaie@uwaterloo.ca \\
University of Waterloo, Waterloo, ON, Canada
}
	
\editor{ }

\maketitle

\begin{abstract}
We propose the concepts of philomatics and psychomatics as hybrid combinations of philosophy and psychology with mathematics. We explain four motivations for this combination which are fulfilling the desire of analytical philosophy, proposing science of philosophy, justifying mathematical algorithms by philosophy, and abstraction in both philosophy and mathematics. We enumerate various examples for philomatics and psychomatics, some of which are explained in more depth. The first example is the analysis of relation between the context principle, semantic holism, and the usage theory of meaning with the attention mechanism in mathematics. The other example is on the relations of Plato's theory of forms in philosophy with the holographic principle in string theory, object-oriented programming, and machine learning. Finally, the relation between Wittgenstein's family resemblance and clustering in mathematics is explained. This paper opens the door of research for combining philosophy and psychology with mathematics. 
\end{abstract}
\hfill\break


\medskip

\section{Introduction and Definition}\label{section_introduction}

In this paper, we propose the concepts of philomatics and psychomatics which are the hybrid combination of philosophy and psychology with mathematics. The formal definitions of these newly fields of research are provided in the following.  

\begin{definition}[Philomatics]
Philomatics is the combination of philosophy and mathematics, in which either the philosophical concepts are expressed and described using mathematical concepts and equations or the mathematical concepts are analyzed philosophically. 
\end{definition}

\begin{definition}[Psychomatics]
Psychomatics is the combination of psychology (including, but not limited to, existential psychology/psychotherapy \cite{may1961existential,yalom1980existential} and psychoanalysis \cite{bateman2021introduction}) and mathematics, in which either the psychological concepts are expressed and described using mathematical concepts and equations or the mathematical concepts are analyzed psychologically. 
\end{definition}

The reason that we separate philomatics and psychomatics as two different proposed fields of study is that, according to Gottlob Frege's criticism against psychologism, philosophy is different from psychology and these two should not be confused with each other \cite{dabbagh2017analytical}.

It is noteworthy that, as mathematics includes various sub-fields, combining philosophy and psychology with any sub-field of mathematics can be considered as philomatics and psychomatics, respectively.
For example, as machine learning is nothing but optimization, linear algebra, probability, and statistics, combinations of philosophy and psychology with machine learning are special categories of philomatics and psychomatics, respectively. Another example is combining philosophical and psychological concepts with the concepts of physics, whose language is mathematics; an example of this will be provided in Section \ref{section_plato_theory_ideas}.

Some brief examples of philomatics and psychomatics are as follows. 
Some philosophies can be stated as optimization problems \cite{boyd2004convex}; this is because some quality needs to be minimized or maximized in philosophy or psychology. If the quality can be quantified, an optimization problem can be formed ans solved using optimization techniques \cite{ghojogh2021kkt}. 
Another possible example is modeling concepts of philosophy or psychology as games in game theory \cite{rasmusen2006games} where the agents either compete or collaborate. The agents can be philosophical or psychological concepts. 
Modeling these concepts as games makes sense because trade-offs can be often observed in real-life concepts which are suitable to be modeled as games. 
For the same reason, multi-agent reinforcement learning \cite{sutton2018reinforcement,bucsoniu2010multi} can be used for philomatics and psychomatics.

As one of the main concentrations of the analytical philosophy is on the philosophy of language \cite{dummett1985origins,glock2008what}, natural language processing \cite{jurafsky2008speech} in machine learning and mathematics can be used for modeling philosophy as another example of philomatics. 
One other possible example of philomatics and psychomatics is the usage of chaos theory \cite{thompson1990nonlinear} to model human's or society's behaviors using partial differential equations as nonlinear dynamical or chaotic systems with different boundary conditions \cite{loye1987chaos,robertson2014chaos}. 

Moreover, as modern physics \cite{tipler2012modern}, including relativity and quantum physics, is highly founded on mathematics, it can be used to model philosophy and psychology. Quantum mind, also called the quantum physics of consciousness, is an example which has tried to model consciousness in the philosophy of mind, by quantum physics. It was first proposed by Eugene Wigner who modeled consciousness as the quantum wave function collapses \cite{wigner1995remarks}. Later, many works followed up his work for quantum mind for unifying physical and social ontology \cite{wendt2015quantum}.

\section{Motivations}

We have had multiple motivations for proposing philomatics and psychomatics, some of which are introduced in the following. 

\subsection{The First Motivation}

It used to be the very desire and wish of analytical philosophy to have a scientific approach to philosophy. 
Initially, this goal was set by Bertrand Russell, Gottlob Frege, and George Edward Moore followed by Ludwig Wittgenstein afterwards. 
Russell and Frege started condensed research on making a universal logical language for better philosophizing \cite{sullivan2003logicism}. In his seminal paper \cite{russell1905denoting}, Russel argues that the current languages that humans use for communication do not suffice for precise philosophization. As an example, he mentions the sentence ``The present king of France is not bold". This sentence seems to have a correct structure but it is meaningless because France does not have any king at the present time.
Later, inspired by this, the early Wittgenstein justified in his Tractatus Logico-Philosophicus \cite{wittgenstein1921tractatus} that the usual language cannot carry the amount of responsibility that the philosophers have expected. 

The Vienna Circle, specifically, had the wish of having a scientific view to philosophy.
The members of the Vienna Circle were philosophers who used to gather in Vienna of Austria and were especially inspired by Wittgenstein's Tractatus Logico-Philosophicus \cite{wittgenstein1921tractatus}. 
In 1922, they announced their manifesto named ``The Scientific World-Conception"\footnote{This manifesto was in German and its title has been translated to various names, such as ``The Scientific World-Conception", ``The Scientific Conception of the World", and ``The Scientific Worldview"} \cite{vienna1922Scientific} stating that metaphysics, discussed in the continental philosophy, is non-sense and there should be a scientific view to the world and philosophy. 

Russell, Frege, and the early Wittgenstein tried to make this desire true by combining philosophy and logic. They tried to have a logical view to philosophy so they purify the language of philosophy from non-sense statements. However, Russell confessed later that the project failed finally and could not succeed. 
By proposing the idea of philomatics, we try to fulfill this old desire in another way by combining philosophy and mathematics, rather than combining philosophy and logic. Although, the mathematical statements have a logical structure, however, mathematics goes beyond logic and this may be the success key for the puzzle of purifying philosophy from non-sense statements. 

\subsection{The Second Motivation}

There are two main approaches in analytical philosophy toward the relation of philosophy and science \cite{dabbagh2017analytical}. One approach states that there is not essential difference between philosophy and science and that philosophy is a a part of science. In this perspective, the task of philosophy is to address the questions which science cannot have an answer for. Bertrand Russell and Willard Van Orman Quine had this perspective. 
The second approach for the relation of science and philosophy believes a qualitative difference between philosophy and science. It says that philosophy is a second-stage activity or a logical science. A few example philosophers having this approach are the early Wittgenstein, the Vienna Circle, the Cambridge school of analysis, the Oxford school of analysis, and Rudolf Carnap.

A dialectic synthesis of these two approaches to the relation of philosophy and science was the development of philosophy of science \cite{bird1998philosophy,burtt2003metaphysical}. The field of philosophy of science has a philosophical perspective to science and explains what qualifies to be a scientific approach. A special case of philosophy of science is the philosophy of mathematics \cite{shapiro1997philosophy}.
By proposing philomatics and psychomatics, we develop another synthesis as a new field which is the other way around, i.e., science of philosophy rather than philosophy of science! As a special case, we propose to develop the mathematics of philosophy rather than philosophy of mathematics.

\subsection{The Third Motivation}

Recently, many methods are being developed in mathematics (and in machine learning as a specific field of mathematics) which are completely based on the previous philosophical discussions but they are not consciously aware of those discussions. It is valuable to inform the society that some algorithms and methods which have been developed were already discussed before in philosophy!
An example of this is the attention mechanism which was proposed in 2017 for natural language processing \cite{vaswani2017attention}. The attention mechanism, which is based inner product for measuring similarity between words in the context, is basically using the ``context principle" which was proposed by Frege in 1884--1892 \cite{frege1884foundations,frege1892sense} more than a century ago!
We will discuss this further in Section \ref{section_attention_mechanism}.

\subsection{The Fourth Motivation}


The fourth motivation for proposing philomatics and psychomatics is that both philosophy and mathematics are abstract. In Frege's philosophy, the meaning and concepts are abstract entities in an abstract world \cite{dabbagh2017analytical}. This may be inspired by Plato's Theory of Ideas \cite{demos1939philosophy} in which ideal entities exist in an abstract world. 
Moreover, the early Wittgenstein believed in the ``picture theory of meaning" in which the meaning of words in a language exist in an abstract world of logic \cite{wittgenstein1921tractatus}. 
These examples show that abstract concepts have been appeared in various types of philosophy, from ancient philosophy to modern analytical philosophy. On the other hand, mathematics also deals with abstract concepts \cite{ferrari2003abstraction}.
This explains that philosophy and mathematics can be related in the sense of abstraction. 

\section{A Few Examples for Philomatics and Psychomatics}

Some examples of possible research in philomatics and psychomatics were enumerated in Section \ref{section_introduction}. Here, we briefly introduce three examples of philomatics in more details. Later in the follow-up papers, we will thoroughly explain these examples in more depth.
The examples in Sections \ref{section_attention_mechanism} and \ref{section_plato_theory_ideas} are sample analytical/comparative works in philomatics while the example in Section \ref{section_family_resemblance} is a sample compositional work which designs a new mathematical algorithm using philosophical concepts.

\subsection{Relation of Context Principle and Semantic Holism with Attention Mechanism}\label{section_attention_mechanism}


\subsubsection{Attention Mechanism and Context-Aware Embedding}

Attention mechanism was primarily noticed in cognitive science at the end of 20-th century \cite{tomlin1994attention}; however, it was mainly proposed for natural language processing recently. 
Before the proposal of the attention mechanism in 2017 \cite{vaswani2017attention}, recurrent neural networks \cite{rumelhart1986learning} and long short-term memory networks \cite{hochreiter1997long} were mostly used for natural language processing. However, attention mechanism changed the research path by using the concept of attention which is a mathematical model for feature weighting. Assume you are looking at the portrait of Mona Lisa. You do not attend to all regions if the picture; rather, you focus on her face more than the river in her background. Likewise, in skimming as a technique for fast reading, you do not attend to all words equally but you focus on important words and ignore the less informative words. This more and less attention to different features of picture or text can be modeled mathematically as weighting \cite{ghojogh2020attention}. 

In the attention mechanism \cite{vaswani2017attention}, the different words of a sentence or a text corpus attend to each other with different weights, where the weights are obtained by the amount of similarity of every word with another word. As a result, the more similar two words are, the more attention they have to each other. Two language models, named BERT (Bidirectional Encoder Representations from Transformers) \cite{devlin2018bert} and GPT (Generative Pre-trained Transformer) \cite{radford2018improving} have been proposed based on the attention mechanism. 
In contrast to the previous language models, such as word2vec \cite{mikolov2013efficient,mikolov2013distributed} and GloVe \cite{pennington2014glove}, BERT and GPT are based context-aware language models which provide different embedding vectors for a word in different contexts. For example, the word ``bank" has different meanings and therefore different embeddings in the sentences ``Money is in the bank" and ``Some plants grow in bank of rivers".

\subsubsection{Relation with the Philosophy of Language}

The explained attention mechanism is closely related to Gottlob Frege's ``context principle" \cite{frege1884foundations,frege1892sense} and Willard Quine's ``semantic holism" \cite{quine1951two} in analytical philosophy and the philosophy of language. Note that the attention mechanism was proposed for natural language processing in 2017 while the context principle was proposed more than a century ago in 1884--1892.
In 1884, Frege stated that ``never [...] ask for the meaning of a word in isolation, but only in the context of a proposition" \cite{frege1884foundations}. Moreover, in 1892, Frege stated that a word acquires its meaning only within the context of an entire sentence \cite{frege1892sense}.  

The context principle or the contextualism \cite{preyer2005contextualism} states that the meaning of a word is understandable only as a part of a whole, i.e., a sentence or context. 
This was developed further into the ``usage theory of meaning" in the philosophy of later Wittgenstein. In his ``Philosophical Investigations", Wittgenstein stated ``comprehending a proposition means comprehending a language" \cite{wittgenstein1953philosophical}.
This ``usage theory of meaning" replaced the ``picture theory of meaning" in his early philosophy which claimed that the meaning of words in a language exist in an abstract world of logic (cf. Wittgenstein's Tractatus Logico-Philosophicus \cite{wittgenstein1921tractatus}).
After the later Wittgenstein, this theory was further developed and generalized into Quine's ``semantic holism" stating ``the unit of measure of empirical meaning is all of science in its globality" \cite{quine1951two}.
In Quine's perspective, the meaning of a claim is perceivable only if it is considered in a sequence of propositions and contexts. In this perspective, human's mind is like a court which considers all the evidences for understanding the meaning of something \cite{dabbagh2017analytical}. 

Comparing the above explanations in natural language processing and the philosophy of language shows that they have close relation. Attention mechanism, as a weighting scheme in mathematics, is using the ideas of context principle, usage theory of meaning, and semantic holism in the philosophy of language. All of these theories are defining the meaning of every word in the natural language based on the context it appears in. 

\subsection{Relation of Plato's Theory of Forms with Concepts in Mathematical Physics, Computer Programming, and Machine Learning}\label{section_plato_theory_ideas}


Another comparative-analytical work in philomatics can be explaining the relation of Plato's theory of forms in philosophy with concepts in mathematical physics, computer programming, and machine learning. In the following, these relations are explained. 

\subsubsection{Plato's Theory of Forms}

In metaphysics, "particulars" are the existing objects in the world and "universals" are the properties, such as colors, which are shared between many particulars \cite{loux2017metaphysics}. 
There are two approaches regarding the ontology (existence) of universals in the world. Realism states that the universals exist in the world while nominalism claims that there are no existing universals and they are nothing but resembled characteristics between the particulars \cite{armstrong1978nominalism}. 

One of the first approaches in realism was proposed by Plato's theory of forms, also called the theory of ideas. 
According to this theory, everything in this physical world has some true essence, or idea (also called form), in a world of ideas where the things merely imitate the true ideas. In this theory, the physical world is not as true as the timeless, absolute, unchangeable ideas \cite{sedley2016introduction}.
In his Republics \cite{plato2008republic}, Plato has illustrated his theory by an example where several people are in a cave facing the wall of the cave. The true ideas are outside of the cave in the real world but their shadows appear on the wall. The people in the cave suppose that the shadows are the real things; however, they are just shadows of the real things existing outside of the cave. 
Some researchers believe that the theory of ideas has had significant impact on the monotheistic religions \cite{elkaisy2009afterlife}. 

The theory of forms is also one of the origins of idealism in metaphysics stating that the reality is entirely a mental construct and only the ideas are real or they are the highest level of reality \cite{goldschmidt2017idealism}. Many philosophers developed other variants of idealism; e.g., Berkeley, Kant, and Hegel developed empirical \cite{berkeley1710treatise}, transcendental \cite{kant1781critique}, and absolute \cite{hegel1807phenomenology} idealism, respectively \cite{dunham2014idealism}. 

\subsubsection{Relation of Plato's Theory of Forms with the Holographic Principle in String Theory}

The theory of forms in philosophy is completely related to the holographic principle in modern physics and string theory. 
Black holes are the regions in the space-time of Einstein's field equations where the gravity, i.e., the curvature of space-time, is too huge that nothing including light cannot escape from. 
This causes an information paradox because the quantum information of the things falling in black hole should not be lost. One possible solution to this paradox is the evaporation of black holes through the Hawking's radiation \cite{hawking1974black}.
Another possible way to resolve this paradox is the holographic principle \cite{stephens1994black} which states that the information of falling things into the black hole is encoded as a lower-dimensional copy on the boundary of the black hole. In this way, the information is not lost but its low-dimensional copy is restored. It is noteworthy that this low-dimensional encoding is also related to dimensionality reduction in mathematical machine learning \cite{ghojogh2023elements}. 
Physicists have even hypothesized that the current world might just be the encoded hologram of another high-dimensional world \cite{susskind1995world}; this makes sense because the world has been assumed to be of higher dimensions than four (three dimensions for space and one dimension for time) in the string theory or the M theory \cite{rizzo2004pedagogical}.

The main idea behind the holographic principle is very similar to the theory of forms by Plato. 
In the holographic principle, the world is supposed to the be the hologram of a high dimensional world. In Plato's theory of forms, also, the world is a projection or shadow of a true world of ideas (forms). 

\subsubsection{Relation of Plato's Theory of Forms with Object Oriented Programming}

There is also a relation between Plato's theory of forms and object oriented programming \cite{barnes2000object} in computer programming. 
In object oriented programming, which is a modular programming suitable for team-working, the programmer defines some classes as modules of code which can perform some tasks. 
The computer program can make objects or instances from these classes. 
Every class has several attributes which are initially set for every individual object  when the object is created. Moreover, every class has some functions for performing some tasks. These functions are also called methods. One of the methods, named constructor, is for setting the initial values of attributes. 
The classes can inherit from each other in the way that the inheritor inherits some methods (functions) and attributes from the other class. 

As it is obvious, there is a close relation of object oriented programming and the theory of forms. The classes are like the forms or ideas in the world of forms. The objects or instances of classes are like the instances of forms in this phenomenal world. Every class or the form has some attributes and the instances of the class or the form has those attributes consequently. 
The classes or the forms can also inherit from each other. 

An example is the form/class of human. Benyamin and Morteza are two instances/objects of human. The form/class of human has some attributes such as height, weight, color, language, mood, etc. Benyamin and Morteza each have some states/values of these attributes. 
These attributes are set initially at the time of birth or by the constructor method, when the instance/object is created. However, these attributes might be reset again later to other values. 
The form/class of human has some functions or methods such as walking, sitting, talking, thinking, etc. 
At every moment, Benyamin and Morteza perform one or several of these functions. 
The form/class of human inherits from the form/class of mammal, which itself inherts from the form/class of animal.
Some functions (methods) and attributes of the form/class of human are inherited from the form/class of mammal. However, some functions and attributes are specific to the form/class of human.
For example, the way of intercourse and giving birth is inherited from mammals but talking in complicated languages and thinking to its own death can be enumerated as the functions specific to human. 

\subsubsection{Relation of Plato's Theory of Forms with Machine Learning}

The theory of forms is also related to classification and clustering in machine learning. 
Consider a set of data points, called a dataset. They can be grouped into multiple clusters or classes of data, depending on whether they have labels or not.  
In both cases of clustering and classification, every cluster or class can be considered as a group. Every group can have a representative. This representative is usually a summary statistics, such as mean, of the group. If the summary statistics is not median, the representative of a group may not be any of the data points in that group. As a result, this representative can be considered as the true form of the group with instances existing in the group. 

Another relation of the theory of forms with machine learning can be considered in generative models. 
In many of the generative models, a latent noise is fed to the model and the model generates some semantic data instance corresponding to that noise \cite{harshvardhan2020comprehensive}. 
The most fundamental generative model in machine learning is factor analysis \cite{fruchter1954introduction,ghojogh2021factor} and one of the state-of-the-art generative models is generative adversarial networks \cite{goodfellow2014generative,ghojogh2021generative}.
In generative models, every latent noise corresponds to a generated data point. As a result, the latent noise can be considered as the true form of the generated instance in the world of forms.
It is noteworthy that vector arithmetics can be performed between the latent noise vectors \cite{radford2016unsupervised}. For example, if the corresponding noise to a facial image without eyeglasses is subtracted from the corresponding noise to a facial image with eyeglasses, the obtained noise vector corresponds to the image of eyeglasses! The corresponding noise to the eyeglasses can be considered as the form of eyeglasses in the world of ideas. 

\subsection{Relation of Wittgenstein's Family Resemblance with Clustering}\label{section_family_resemblance}

A compositional example for philomatics is designing a clustering algorithm based on Wittgenstein's ``family resemblance". 
In his ``Philosophical Investigations", the later Wittgenstein proposed the concept of family resemblance explained in the following. Consider three entities 1, 2, and 3. Suppose the entity 1 has attributes a and b, the entity 2 has attributes b and c, and the entity 3 has attributes c and d. While the first and third entities do not have any attribute in common, all the three entities can be considered to belong to a family. This is because entities 1 and 2 have something in common and entities 2 and 3 have something in common. As a result, they have resemblance in a chain or tree structure. 

A careful mathematical look at Wittgenstein's ``family resemblance" shows that it is a special algorithm for clustering \cite{xu2005survey} in which the members of every cluster have similarity or resemblance in a chain structure. As this similarity is in the form of a chain or tree, graph theory \cite{west2001introduction} can also be employed to design such a clustering algorithm. 
Consider several vectors in a Euclidean space.
Any measure of similarity, such as inner product, can be used to measure the resemblance of vectors. Using a threshold, it can be determined if every two vectors are resembled or not. i.e., if they have any common attribute or not. Then, using the idea of chain-resemblance in Wittgenstein's family resemblance, the clusters can be detected. In conclusion, Wittgenstein's family resemblance can be seen as a clustering algorithm in mathematics and a clustering method can be designed according to it. 

\section{Conclusion and Future Directions}

In this paper, we proposed philomatics and psychomatics as two fields of research and study. This paper opens the door of research for combining philosophy and psychology with mathematics. 
In this field of research, various tasks can be performed some of which are mentioned in the following. 
Firstly, the existing theories of philosophy can modeled mathematically by going through the philosophy of previous philosophers, such as Kant, and modeling their philosophy in mathematics. In other words, the language of previous philosophy can be translated to the language of mathematics. 

Secondly, new philosophical and/or psychological theories can be developed using philomatics and psychomatics. For this, the new theories can be explained both in words (i.e., the traditional approach) and in mathematics (i.e., the new approach).
Alternatively, the new philosophical theories can be developed purely with the language of mathematics.
We authors recommend that the combination of words and mathematics can be a suitable approach for developing new philosophical and psychological theories. 

Thirdly, new mathematical algorithms can be developed inspired by the philosophical and psychological theories. An example for this third task is clustering by Wittgenstein's family resemblance (see Section \ref{section_family_resemblance}). Another example is ``affective manifolds" which has been proposed for modeling machine's mind using machine learning based on philosophical and psychological discussions \cite{ghojogh2022affective}. 

\section*{Acknowledgement}

The authors hugely thank Dr. Soroush Dabbagh for the philosophical discussions and his lectures in philosophy.

\bibliography{ref}



\end{document}